\documentclass[11pt]{amsart}
\usepackage{stmaryrd}
\usepackage{mathrsfs}
\usepackage{amssymb}

\usepackage{titletoc}
\input xy
  \xyoption{all}
\usepackage{amscd}
\usepackage{amsmath, amssymb}
\usepackage{amsfonts}
\usepackage[colorlinks,linkcolor=blue,citecolor=blue, pdfstartview=FitH]{hyperref}

\numberwithin{equation}{section}

\theoremstyle{plain}
\newtheorem{thm}{Theorem}[section]
\newtheorem{theorem}[thm]{Theorem}
\newtheorem{lemma}[thm]{Lemma}
\newtheorem{corollary}[thm]{Corollary}
\newtheorem{proposition}[thm]{Proposition}

\theoremstyle{definition}

\newtheorem{question}[thm]{Question}

\newtheorem{remark}[thm]{Remark}

\newtheorem{definition}[thm]{Definition}

\newtheorem{example}[thm]{Example}

\newtheorem{defn-thm}[thm]{Definition-Theorem}
\newtheorem{conjecture}[thm]{Conjecture}

\newcommand{\sO}{{\mathcal O}}

\newcommand{\C}{{\mathbb C}}

\renewcommand{\P}{{\mathbb P}}

\newcommand{\R}{{\mathbb R}}

\newcommand{\qtq}[1]{\quad\mbox{#1}\quad}

\newcommand{\Om}{\Omega}

\newcommand{\ts}{\otimes}

\newcommand{\btheorem}{\begin{theorem}}
\newcommand{\etheorem}{\end{theorem}}
\newcommand{\bproposition}{\begin{proposition}}
\newcommand{\eproposition}{\end{proposition}}
\newcommand{\bdefinition}{\begin{definition}}
\newcommand{\edefinition}{\end{definition}}
\newcommand{\bcorollary}{\begin{corollary}}
\newcommand{\ecorollary}{\end{corollary}}
\newcommand{\bproof}{\begin{proof}}
\newcommand{\eproof}{\end{proof}}
\newcommand{\bremark}{\begin{remark}}
\newcommand{\eremark}{\end{remark}}
\newcommand{\eexample}{\end{example}}
\newcommand{\bexample}{\begin{example}}

\newcommand{\elemma}{\end{lemma}}
\newcommand{\blemma}{\begin{lemma}}

\renewcommand{\bar}{\overline}
\newcommand{\eps}{\varepsilon}

\renewcommand{\phi}{\varphi}

\newcommand{\ee}{\end{eqnarray*}}
\newcommand{\be}{\begin{eqnarray*}}

\newcommand{\beq}{\begin{equation}}
\newcommand{\eeq}{\end{equation}}

\newcommand{\bd}{\begin{enumerate}}
\newcommand{\ed}{\end{enumerate}}

\renewcommand{\tilde}{\widetilde}

\renewcommand{\>}{\rightarrow}

\usepackage{fancyhdr}
\renewcommand{\leftmark}{{\scriptsize K\"ahler manifolds homotopic to negatively curved  manifolds }}

\begin{document}
\title{Compact K\"ahler manifolds homotopic to negatively curved Riemannian manifolds }
\makeatletter
\let\uppercasenonmath\@gobble
\let\MakeUppercase\relax
\let\scshape\relax
\makeatother
\author{Bing-Long Chen and Xiaokui Yang}

\address{{Address of Bing-Long Chen: Department of Mathematics,
Sun Yat-sen University, Guangzhou, 510275, China.}}
\email{\href{mailto:mcscbl@mail.sysu.edu.cn}{{mcscbl@mail.sysu.edu.cn}}}
\address{{Address of Xiaokui Yang: Morningside Center of Mathematics, Institute of
        Mathematics, Hua Loo-Keng Key Laboratory of Mathematics,
        Academy of Mathematics and Systems Science,
        Chinese Academy of Sciences, Beijing, 100190, China.}}
\email{\href{mailto:xkyang@amss.ac.cn}{{xkyang@amss.ac.cn}}}
\maketitle

\begin{abstract} In this paper, we show that any compact K\"ahler
manifold homotopic to a compact Riemannian manifold with negative
sectional curvature  admits  a K\"ahler-Einstein metric  of general
type. Moreover, we prove that,  on a compact symplectic manifold $X$
homotopic to a compact Riemannian manifold with negative sectional
curvature, for any almost complex structure $J$ compatible with the
symplectic form, there is  no non-constant $J$-holomorphic entire curve $f:\C\> X$.
\end{abstract}

\section{Introduction}

In 1970s, S.-T. Yau proposed the following conjecture:

\vskip 0.3\baselineskip

\noindent\textbf{Conjecture.}  Let $(X,\omega)$ be a compact
K\"ahler manifold with $\dim_\C X>1$. Suppose $(X,\omega)$ has
negative Riemannian sectional curvature, then $X$ is rigid, i.e. $X$
has only one  complex structure.

\vskip 0.5\baselineskip


\noindent It is a fundamental problem on the rigidity of K\"ahler
manifolds with negative curvature.  Yau proved in
\cite[Theorem~6]{Yau77} that when $X$ is covered by a $2$-ball, then
any complex surface oriented homotopic to $X$ must be biholomorphic to $X$. By using the terminology of ``strongly
negativity", Siu established in \cite[Theorem~2]{Siu80} that a
compact
 K\"ahler manifold of the same homotopy type as  a compact K\"ahler manifold
$(X,\omega)$ with strongly negative curvature and $\dim_{\C}X>1$
must be either biholomorphic or conjugate biholomorphic to $X$. Note
that when $\dim_{\C}X=2$, the above Yau's conjecture has been
completely solved by Zheng \cite{Zhe95}.  It is well-known that the
strongly negative curvature condition can imply the negativity of
the Riemannian sectional curvature.

\vskip 0.3\baselineskip

Based on Yau's conjecture, one can also ask a more general question:
if a K\"ahler manifold--or complex manifold--$X$ admits a Riemannian
metric with negative sectional curvature, is there any restriction
on the complex structure of $X$?

\vskip 0.3\baselineskip

The  first main result of our paper is one important step towards the question:

\btheorem\label{main}  Let $X$ be a compact  manifold homotopic to a
compact Riemannian manifold $Y$ with negative sectional curvature.
If $X$ has a K\"ahler complex structure $(J,\omega)$,  then it
admits a K\"ahler-Einstein metric of general type. Moreover, each
submanifold of $(X,J)$  admits a K\"ahler-Einstein metric of general
type. \etheorem

\noindent Indeed, we prove in Theorem \ref{general} a more general
result. One of the main ingredients in the proofs is a notion called
``K\"ahler hyperbolicity" introduced by Gromov in  \cite{Gro91}, and
our key observation is that every K\"ahler hyperbolic manifold has
ample canonical bundle (see Theorem \ref{GH}), which also answers a question asked by
Gromov in \cite[p.267]{Gro91}.

\vskip 0.5\baselineskip

Recall that a  compact  complex manifold $X$ is called Kobayashi (or
Brody) {hyperbolic} if every holomorphic map $f:\C\>X$ is constant.
Gromov pointed out in \cite{Gro91} that  every K\"ahler hyperbolic
manifold is Kobayashi hyperbolic.  On the other hand, one can  also  extend these terminologies to symplectic manifolds. Note that for a fixed symplectic form $\omega$, there are many almost complex structures compatible with $\omega$.  Our second result is

\btheorem\label{symp}  Let $(X,\omega)$ be a compact symplectic
manifold homotopic to a compact Riemannian manifold with negative
sectional curvature. For any almost complex structure $J$ on $X$
compatible with $\omega$, there exists no non-constant
$J$-holomorphic map $f: \C\>X$. \etheorem

\noindent \textbf{Acknowledgements.} B.-L. Chen was partially
supported by grants NSFC11521101, 11025107.  X.-K. Yang was
partially supported by China's Recruitment
 Program of Global Experts and National Center for Mathematics and Interdisciplinary Sciences,
 Chinese Academy of Sciences. The authors would like to thank the anonymous
 referee for  pointing out a similar proof in Lemma
 \ref{sectional}.

\section{Background materials}

\subsection{Positivity of line bundles}

 Let $(X,\omega)$ be a  smooth projective manifold of complex dimension $n$, $L\>X$ a
holomorphic line bundle and $E\>X$ a holomorphic vector bundle. Let
$\sO_{\P(E^*)}(1)$ be the tautological line bundle of the projective
bundle $\P(E^*)$ over $X$. \bd\item $L$ is said to be
 \emph{ample}  if $L^k$ is
very ample for some large $k$, i.e. the map $X\>\P(H^0(X,L^k)^*)$
defined by the global sections of $L^k$  is a holomorphic embedding.
$L$ is called \emph{semi-ample} if for some
large positive integer $k$, $L^{ k}$ is  generated by its global
sections, i.e. the evaluation map $\iota: H^0(X,L^{ k})\>L^{ k}$ is
surjective. The vector bundle $E$ is called \emph{ample}   if $\sO_{\P(E^*)}(1)$
is an ample  line bundle.

\item $L$ is said to be \emph{nef} ( or \emph{numerically effective}), if $L\cdot C\geq
0$ for any irreducible curve $C$ in $X$. 

\item $L$ is said to be \emph{big}, if the Kodaira dimension $\kappa(L)=\dim X$ where $$\kappa(L):=\limsup_{m\>+\infty} \frac{\log \dim_\C
H^0(X,L^{ m})}{\log m}.$$  Here we use the convention that the
logarithm of zero is  $-\infty$.
 \ed

\bdefinition  $X$ is said to be \emph{of general type} if the
Kodaira dimension $\kappa(X):=\kappa(K_X)$ is equal to the complex
dimension of $X$. \edefinition

\noindent  There are many examples of compact complex manifolds of
general type. For instance, manifolds with ample canonical bundles.

\subsection{K\"ahler hyperbolicity} Let's recall some concepts
introduced by Gromov in \cite{Gro91}. Let $(X,g)$ be a Riemannian
manifold. A differential form $\alpha$ is called \emph{$d$-bounded}
if there exists a  form $\beta$ on $X$ such that $\alpha=d\beta$ and
\beq \label{2.1} \|\beta\|_{L^\infty(X,g)}=\sup_{x\in X} |\beta(x)|_g<\infty. \eeq
It is obvious that if $X$ is compact, then every exact form is
$d$-bounded. However, when $X$ is not compact, there exist smooth
differential forms which are exact but not $d$-bounded. For
instance, on $\R^n$,  $\alpha=dx^1\wedge\cdots \wedge d x^n$ is
exact, but it is not $d$-bounded.

\bdefinition Let $(X,g)$ be a Riemannian manifold and $\pi:(\tilde
X, \tilde g )\>(X,g)$ be the universal covering with $\tilde g=\pi^*
g$. A form $\alpha$ on $X$ is called \emph{$\tilde d$-bounded} if
$\pi^*\alpha$ is a $d$-bounded form on $(\tilde X, \tilde g)$.
\edefinition

\noindent It is obvious that the $\tilde d$-boundedness does not
depend on the metric $g$ when X is compact.

 \blemma Let $(X,g)$ be a compact
Riemannian manifold. If $\alpha$ is $\tilde d$-bounded on $(X,g)$, then for any
metric $g_1$ on $X$, $\alpha$ is also $\tilde d$-bounded on
$(X, g_1)$. \elemma

\bproof Since $X$ is compact, any two smooth metrics on $X$ are
equivalent. \eproof

\noindent Note also that the $\tilde d$-boundedness of a closed form
$\alpha$ on a compact manifold $X$ depends only on the cohomology class $[\alpha]\in H_{DR}^*
(X,\R)$.

\blemma Let $(X,g)$ be a compact Riemannian manifold. Suppose
$\alpha$ is $\tilde d$-bounded, then $\alpha_1=\alpha+d\gamma$ is
also $\tilde d$-bounded. \elemma

\bproof Let $\pi:(\tilde X, \tilde g)\>(X,g)$ be the universal
covering and $\beta$ be the form on $\tilde X$ such that
$\pi^*\alpha=d\beta$ and $\|\beta\|_{L^\infty(\tilde X, \tilde
g)}<\infty$. Hence, we have
$$\pi^*\alpha_1=d(\beta+\pi^*\gamma)$$
and \be \|\beta+\pi^*\gamma\|_{L^\infty(\tilde X, \tilde g)}&\leq&
\|\beta\|_{L^\infty(\tilde X, \tilde g)}+
\|\pi^*\gamma\|_{L^\infty(\tilde X, \tilde g)}\\
&= & \|\beta\|_{L^\infty(\tilde X, \tilde g)}+
\|\gamma\|_{L^\infty( X,  g)}<\infty. \ee
 \eproof


\bdefinition Let $X$ be a  Riemannian manifold.  $X$ has
\emph{$\tilde d$-bounded $i^{\text{th}}$ cohomology} if every class
in $H_{\text{DR}}^i(X,\R)$ is $\tilde d$-bounded. \edefinition

\blemma\label{pullback} Let $f:X\>Y$ be a smooth map between two
compact Riemannian manifolds. Suppose $\alpha$ is $\tilde d$-bounded
on $Y$, then $f^*\alpha$ is $\tilde d$-bounded on $X$. \elemma

\bproof  Let $\pi_X:\tilde X\>X$ and $\pi_Y:\tilde Y\>Y$ be the
universal coverings of $X$ and $Y$ respectively. Since $\tilde X$ is
simply connected, there exists a lifting map $\tilde f:\tilde
X\>\tilde Y$, such that the following diagram
$$ \xymatrix{
  \tilde X \ar[d]_{\pi_X} \ar[r]^{\tilde f}
                & \tilde Y \ar[d]^{\pi_Y}  \\
  X\ar[r]_{f}
                & Y            }$$
 commutes.  On the other hand, we know $\pi_Y^*\alpha=d\beta$
for some $L^\infty$-bounded form $\beta$ over $\tilde Y$. Hence \beq
\pi_X^*(f^*\alpha)=\tilde f^*(\pi_Y^*\alpha)=\tilde
f^*(d\beta)=d(\tilde f^*\beta).\eeq Since $X$ and $Y$ are compact,
$\pi_X$ and $\pi_Y$ are local isometries, \beq \|\tilde
f^*\beta\|_{L^\infty(\tilde X, \pi^*_X g_X)}\leq C
\|\beta\|_{L^\infty(\tilde Y, \pi^*_Y g_Y)}\cdot
\|f\|^{p}_{C^1(g_X,g_Y)}<\infty\label{sequence2}\eeq where $p$ is the degree of $\beta$ and $C$ is a
constant depending only on $X$ and $p$.
 \eproof

In geometry, various notions of hyperbolicity have been introduced,
and the typical examples are manifolds with negative curvature in
suitable sense. The starting point for the present investigation is
Gromov's notion of K\"ahler hyperbolicity \cite{Gro91}.

\bdefinition\label{KH} Let $X$ be a compact complex manifold. $X$ is
called \emph{K\"ahler hyperbolic} if it admits a K\"ahler metric
$\omega$ such that $\omega$ is $\tilde d$-bounded. \edefinition

\noindent The typical  examples of K\"ahler hyperbolic manifolds are
locally Hermitian symmetric spaces of noncompact type.

\vskip 0.4\baselineskip

\noindent As we mentioned before, a compact complex manifold is
called Kobayashi hyperbolic if it contains no entire curves.    A fundamental problem
in complex  geometry is Kobayashi's  conjecture (e.g. Lang's survey paper
\cite{Lan86}):
\begin{conjecture} Let $X$ be a compact complex manifold. If $X$ is
  Kobayashi hyperbolic, then the canonical bundle $K_X$ is ample.
 \end{conjecture}
  \noindent Along the
same line, Gromov asked the following question in
\cite[p.267]{Gro91}.

\begin{question} Let $X$ be a compact  K\"ahler hyperbolic
manifold. Is the canonical bundle $K_X$ ample? Is the cotangent
bundle $\Om_X^1$ ample?
\end{question}

\noindent We first give a counter-example to the second part of
Gromov's question. \bexample Let $X=C_1\times C_2$ be the product of
two smooth curves of genus at least $2$. It is obvious that $X$ is
K\"ahler hyperbolic since both $C_1$ and $C_2$ are K\"ahler
hyperbolic. The cotangent bundle is $\Om_X^1=\pi^*_1\Om^1_{C_1}\ts
\pi_2^*\Om^1_{C_2}$, which is not ample. Indeed, its restriction to
a curve $C_1\times \{p\}$ has a trivial summand. \eexample

\noindent Next, we give an affirmative answer to the first part of
Gromov's question based on an observation in algebraic geometry. To
the readers' convenience, we include a straightforward proof here.

\btheorem\label{GH} Let $X$ be a compact   K\"ahler hyperbolic manifold. Then
the canonical bundle $K_X$ is ample.

\etheorem \bproof If $X$ is K\"ahler hyperbolic, then
  $X$ contains no rational curves. Indeed, suppose $f:\P^1\>X$ is a rational curve. We want to show $f$ is
a constant, i.e. $f^*\omega=0$.   Let $\pi_X:\tilde X\>X$ be the
universal covering. Then there is a lifting $\tilde f:\P^1\>\tilde
X$ such that $\pi_X\circ \tilde f =f$. Since $\omega$ is $\tilde
d$-bounded, i.e. there exists a bounded $1$-form $\beta$ on $\tilde
X$ such that $\pi^*_X\omega=d\beta$,  \beq f^*\omega=\tilde
f^*(\pi^*_X \omega)= d(\tilde f^*\beta).\eeq It implies
$$\int_{\P^1}f^*\omega=\int_{\P^1} d(\tilde f^*\beta)=0,$$
and so $f^*\omega=0$.

Gromov proved in \cite[Corollary ~0.4C]{Gro91} that if $X$ is
K\"ahler hyperbolic, then $K_X$ is a big line bundle, and so $X$ is
 Moishezon. By Moishezon's theorem, the K\"ahler and Moishezon manifold $X$  is projective.  Since $X$
 contains no rational curves,  Mori's cone theorem implies that
$K_X$ is nef. Since $K_X$ is big and nef,  by Kawamata-Reid-Shokurov
base point free  theorem,  $K_X$ is semi-ample. Then there exists
$m$ big enough such that $\phi=|m K_X|$ is a morphism. Since $K_X$
is big, there is a positive integer $\tilde m$ such that
$$\tilde mK_X=D+L$$
where $D$ is an effective divisor and $L$ is an ample line bundle.
Suppose $K_X$ is not ample, then there  exists a curve $C$
contracted by $\varphi$, i.e., $K_X\cdot C=0$. Therefore,
$$ D\cdot C=-L\cdot C<0.$$ Let $\Delta=\eps D$ for some small $\eps>0$, then
$(X,\Delta)$ is a klt pair and $K_X+\Delta$ is not $\varphi$-nef.
Then by the relative Cone theorem(e.g. \cite[Theorem~3.25]{KM98})
for log pairs, there exists a rational curve $\tilde{C}$ contracted
by  the morphism $\phi$. This is a contradiction since we have
already proved that $X$ contains no rational curves.  Therefore, we
conclude that  $K_X$ is ample. \eproof

\vskip 1\baselineskip
\section{The proof of Theorem \ref{main}}

In this section, we prove  Theorem \ref{main}, which is based on
 the following result.

\btheorem\label{general} Let $X$ be a compact K\"ahler manifold and
$Y$ be a compact Riemannian manifold with $\tilde d$-bounded
$H_{\emph{DR}}^2(Y,\R)$. Suppose there exist two smooth maps
$f_1:X\>Y$ and $f_2:Y\>X$ such that the image of the induced map
$$(f_2\circ f_1)^{\ast}: H_{\emph{DR}}^2(X,\R)\rightarrow
H_{\emph{DR}}^2(X,\R) $$  contains at least one K\"ahler class. Then $X$ admits
a K\"ahler-Einstein metric  of general type.  Moreover, each
submanifold  of $X$ is also a K\"ahler-Einstein  manifold of general
type.  \etheorem

\bproof Suppose $\tilde X$ and $\tilde Y$ are the universal
coverings of $X$ and $Y$ respectively. Let $\tilde f_1:\tilde
X\>\tilde Y$ and $\tilde f_2:\tilde Y\>\tilde X$ be the liftings of
$f_1$ and $f_2$ respectively such that   the following diagram \beq
\xymatrix{
  \tilde X \ar[d]_{\pi_X} \ar[r]^{\tilde f_1}
                & \tilde Y \ar[d]^{\pi_Y}\ar[r]^{\tilde f_2} &\tilde X\ar[d]^{\pi_X} \\
  X\ar[r]_{f_1} & Y\ar[r]_{f_2}&X } \label{commutative} \eeq
  commutes.
  Let $\omega$ be a K\"ahler metric on $X$ such that $[\omega]$ is   contained in  the image of $$(f_2\circ f_1)^{\ast}: H_{\emph{DR}}^2(X,\R)\rightarrow H_{\emph{DR}}^2(X,\R).$$ Then there
exist a $1$-form $\gamma$ and a closed $2$-form $\omega_1$ on $X$ such
that \beq \label{3.2}\omega=(f_2\circ f_1)^*\omega_1+d\gamma. \eeq Since $Y$
has $\tilde d$-bounded $H_{\text{DR}}^2(Y,\R)$,  for the $2$-form
 $\omega_1$ on $X$, there exists a $1$-form $\beta$ on $\tilde Y$
such that \beq \pi^*_Y\circ f_2^*\omega_1=d\beta\eeq and $\beta$ is
$d$-bounded on $(\tilde Y,\pi^*g_Y)$. It implies $\tilde f_2^*\circ
\pi_X^*\omega_1=d\beta$ and \beq \tilde f_1^*\circ \tilde f_2^*\circ
\pi_X^*\omega_1=\tilde f_1^*d\beta=d(\tilde
f_1^*\beta).\label{sequence} \eeq
Moreover, by (\ref{3.2}) and (\ref{sequence}), we have  \beq
\pi_X^*\omega=d(\pi_X^*\gamma)+\pi^*_X\circ f_1^*\circ
f_2^*\omega_1=d(\pi_X^*\gamma)+\tilde f_1^*\circ \tilde f_2^*\circ
\pi_X^*\omega_1=d(\pi_X^*\gamma+\tilde f_1^*\beta). \eeq By using a similar argument
as in the proof of Lemma \ref{pullback}, we know $\tilde f_1^*\beta$
is bounded on $\tilde X$. Hence, $\pi_X^*\omega$ is $d$-bounded on
$\tilde X$. By definition \ref{KH}, $(X,\omega)$ is K\"ahler hyperbolic. By Theorem \ref{GH}, $K_X$ is ample, i.e. $c_1(X)<0$. Thanks to the
Aubin-Yau theorem, there exists a
 smooth K\"ahler metric $\tilde\omega$ on $X$ such that $Ric(\tilde
 \omega)=-\tilde \omega$.

\vskip 0.4\baselineskip

 Suppose $Z$ is a submanifold of $X$. Let $\omega_Z$ be the K\"ahler metric induced from $(X,\omega)$.
 By Lemma \ref{pullback},
 $(Z,\omega_Z)$ is also K\"ahler hyperbolic. By Theorem \ref{GH} and the Aubin-Yau theorem again,   $Z$ is  a K\"ahler-Einstein manifold of general
 type. \eproof

Before giving the proof of Theorem \ref{main},  we need the
following classical  fact pointed out by Gromov\cite{Gro91} and for
the reader's convenience, we include a simple proof here (see also
\cite[Proposition~8.4]{Bal06}):

\begin{lemma}\label{sectional} Let $(M,g)$ be a simply-connected $n$-dimensional complete Riemannian manifold
with sectional curvature bounded from above by a negative constant,
i.e.  $$\sec \leq- K$$ for some $K>0$. Then for any bounded and closed $p$-form
$\omega$ on $M$, where $p>1$, there exists a bounded $(p-1)$-form
$\beta$ on $M$ such that \beq \omega=d\beta \qtq{and}
|\beta|_{L^{\infty}} \leq {{K}}^{-\frac{1}{2}} |\omega|_{L^{\infty}},\eeq where the
$L^{\infty}$-norm  is given by
\beq \label{3.7} |\beta|_{L^{\infty}}=\sup\left\{ |\beta(v_1,\cdots,v_{p-1})|(x):
x\in M, v_i\in T_{x}M,|v_i|_g=1, i=1,\cdots, p-1\right\}.\eeq
\end{lemma}
(Note that  the norms defined in (\ref{2.1}) and (\ref{3.7}) are
equivalent.)
\begin{proof} Fix $x_0\in M$,  let $\exp_{x_0}: T_{x_0}M\rightarrow M$ be the exponential map,
 which is a diffeomorphism by Cartan-Hadamard theorem. Let $\phi_{t}:M\rightarrow M$, $t\in [0,1]$,
 be a family of maps defined by $\phi_{t}(x)= \exp_{x_0}(t\cdot\exp_{x_0}^{-1}(x))$, $x\in M$. We denote the
 distance function from $x_0$ by $\rho$, then $$X_{t}\mid_{\phi_t(x)}=\left(\frac{d}{dt} \phi_t\right)\mid_{\phi_t(x)}=\rho(x) \nabla \rho\mid_{\phi_t(x)}.$$
  It is clear that $\phi_1=id$ and $\phi_0\equiv x_0$. Then
\begin{equation}
\omega(x)=\int_{0}^{1}
\left(\frac{d}{dt}\phi_{t}^{\ast}\omega\right)(x)
dt=\int_{0}^{1}\phi_t^{\ast} (\mathcal{L}_{X_t}\omega)(x)
dt=d\left(\int_{0}^{1} \phi_{t}^{\ast}(i_{X_t}\omega) dt\right)
\end{equation}
where we have used Cartan's homotopy formula
$\mathcal{L}_{X_t}=d\circ i_{X_t}+i_{X_t}\circ d$
 for differential forms.
If we set \beq \beta=\int_{0}^{1} \phi_{t}^{\ast}(i_{X_t}\omega)
dt,\label{beta}\eeq then $\omega=d\beta$.  We show $\beta$ has
bounded $L^\infty$-norm. Fix $x\in M$, $v_1, v_2,\cdots v_{p-1}\in
T_{x}M$, $|v_i|=1$, $\langle v_i, \nabla \rho\rangle=0$, we have
 \begin{equation}\label{3.10}
 \beta(v_1,\cdots, v_{p-1})(x)=\int_{0}^{1} \omega\left(X_t, (d\phi_{t})(v_1), \cdots, (d\phi_{t})(v_{p-1})\right)(\phi_t(x))dt.
 \end{equation}
 By the standard comparison theorem (e.g. \cite[Theorem~1.28]{CE75}), we have
 \begin{equation}
 |(d\phi_t)(v)|\leq \frac{\sinh (t \sqrt{K}  \rho(x))}{\sinh (\sqrt{K} \rho(x)) }
 \end{equation}
for $v\in T_{x}M$, $|v|=1$ and $\langle v,\nabla \rho\rangle=0$. Hence,
 \begin{equation}
|\beta(v_1,\cdots, v_{p-1})(x)|\leq |\omega|_{L^{\infty}}
\int_{0}^{1} {\rho(x)}\left[\frac{\sinh (t \sqrt{K}  \rho(x))}{\sinh
(\sqrt{K} \rho(x)) }\right]^{p-1} dt.
 \end{equation}
 If $\rho(x)\geq {{K}}^{-\frac{1}{2}}$,
since
 \begin{equation}
 \int_{0}^{1} \sinh^{p-1} (t \sqrt{K}  \rho(x)) dt\leq \frac{\cosh (\sqrt{K}\rho(x))-1}{\sqrt{K}\rho(x)}(\sinh(\sqrt{K}\rho(x)))^{p-2},
 \end{equation}
 we have
 \begin{equation}
|\beta(v_1,\cdots, v_{p-1})(x)|\leq
\frac{|\omega|_{L^{\infty}}}{\sqrt{K}}\cdot\frac{\cosh(\sqrt{K}\rho(x))-1}{\sinh(\sqrt{K}\rho(x))}
\leq {{K}}^{-\frac{1}{2}}|\omega|_{L^{\infty}}.
 \end{equation}
 If $\rho(x)\leq {{K}}^{-\frac{1}{2}}$, we have
 \begin{equation}
\int_{0}^{1} {\rho(x)}\left[\frac{\sinh (t \sqrt{K} \rho(x))}{\sinh
(\sqrt{K} \rho(x)) }\right]^{p-1} dt \leq \rho(x)\leq
{{K}}^{-\frac{1}{2}}.
\end{equation}
Combining two cases, we get $|\beta(v_1,\cdots, v_{p-1})(x)| \leq
{K}^{-\frac{1}{2}} |\omega|_{\infty}$.

  On the other hand, if $v_i$ is parallel to $\nabla \rho$  for some
  $i$, then $(d\phi_{t})(v_i)$ is parallel to $\nabla \rho$. By the explicit  formula (\ref{3.10}),  we see $\beta(v_1,\cdots, v_{p-1})(x)=0$. Hence we obtain
$|\beta|_{L^{\infty}} \leq K^{-\frac{1}{2}} |\omega|_{L^{\infty}}$.
\end{proof}

\vskip 0.5\baselineskip

\noindent \emph{The proof of Theorem \ref{main}.}  By Lemma
\ref{sectional}, we see that a compact Riemannian manifold $Y$ with
negative Riemannian sectional curvature has $\tilde d$-bounded
$q^{\text{th}}$ cohomology for all  $q\geq 2 $.  If $X$ is homotopic
to $Y$, there exist two smooth maps $f_1:X\>Y$ and $f_2:Y\>X$ such
that the induced map $(f_2\circ f_1)^{\ast}$ is identity on
$H_{\emph{DR}}^{\ast}(X,\R)$.  Hence, we can apply Theorem
\ref{general} to complete the proof of Theorem \ref{main}. \qed

\vskip 0.3\baselineskip

As an application of Theorem \ref{main}, we give a slightly shorter
proof on the following rigidity theorem which  was proved by S.-T.
Yau in \cite[Theorem~6]{Yau77}.

\bcorollary[Yau]\label{Yau} Let $N$ be a compact complex surface covered by the unit ball in $\C^2$. Then any complex surface
$M$ that is oriented homotopic to $N$ is biholomorphic to $N$.
\ecorollary \bproof Since $M$ is homotopic to the K\"ahler manifold
$N$ with even first  Betti number $b_1$, $M$ also has even first
Betti number $b_1$  and so it is K\"ahler. Since $N$  has a smooth
metric with strictly negative Riemannian sectional curvature, by
Theorem \ref{main},  $M$ is a K\"ahler-Einstein manifold of general
type. Hence, the classical Chern-number inequality (e.g.
\cite[Theorem~4]{Yau77}) implies ,
$$3c_2(M)\geq c_1^2(M).$$
 Since $M$ is oriented homotopic to $N$,  the signature of $M$ equals that of $N$. One can see
$c^2_1(M) = c^2_1(N)$ and $c_2(M) = c_2(N)$. Therefore, $3c_2(M) =
c_1^2(M)$. By  \cite[Theorem~4]{Yau77}, $M$ is covered by the unit
ball in $\C^2$. By Mostow's rigidity theorem (\cite{Mos73}), $M$ is
in fact biholomorphic to $N$. \eproof

\vskip 1\baselineskip

\section{Hyperbolicity on compact symplectic manifolds}\label{example}

We begin by recalling some basic definitions.  A symplectic form
$\omega$ on a manifold $X$ \emph{tames} an almost complex structure
$J$ if at each point of $X$, $\omega(Z, JZ) >0$ for all nonzero
vectors $Z$. We can define a Riemannian metric by \beq
g_{\omega}(Y,Z)=\frac{1}{2}(\omega(Y,JZ)+\omega(Z,JY)).
\label{smetric}\eeq If, in addition, $$\omega(JY, JZ) =
\omega(Y,Z)$$ for all vectors $Y, Z$, then we say that $\omega$ is
\emph{compatible} with $J$. We establish a more general result than
Theorem \ref{symp}.

\btheorem\label{general2} Let $X$ be a compact symplectic manifold
and $Y$ be a compact Riemannian manifold with $\tilde d$-bounded
$H_{\emph{DR}}^2(Y,\R)$. Suppose there exist two smooth maps
$f_1:X\>Y$ and $f_2:Y\>X$ such that the image of the induced map
$$(f_2\circ f_1)^{\ast}: H_{\emph{DR}}^2(X,\R)\rightarrow
H_{\emph{DR}}^2(X,\R) $$  contains at least one  symplectic class $[\omega]$.
Then for any almost complex structure $J$ on $X$ compatible with the
symplectic form $\omega$,   there exists no non-constant
$J$-holomorphic map $f: \C\ \rightarrow X$. \etheorem

\bproof Let $\omega$ be a symplectic form on $X$ such that
$[\omega]$ is contained in the image of $$(f_2\circ f_1)^{\ast}:
H_{\emph{DR}}^2(X,\R)\rightarrow H_{\emph{DR}}^2(X,\R).$$ Then there
exist a  $1$-form $\gamma$ and a closed $2$-form $\omega_1$ on $X$
such that \beq \omega=(f_2\circ f_1)^*\omega_1+d\gamma. \eeq Now we
use the same commutative diagram as in (\ref{commutative}). There
exists a $1$-form $\beta$ on $\tilde Y$ such that $ \pi^*_Y\circ
f_2^*\omega_1=d\beta$ and $\beta$ is $d$-bounded on $(\tilde
Y,\pi^*g_Y)$. Moreover, we have
$\pi_X^*\omega=d(\pi_X^*\gamma+\tilde f_1^*\beta).$
  If we set $\eta= \pi_X^*\gamma+\tilde
f_1^*\beta$, then $\eta$ is bounded on $\tilde X$, i.e., $\omega$ is
$\tilde{d}$-bounded.

\vskip 0.4\baselineskip

Let $J$ be an almost complex structure on $X$ which is compatible
with the symplectic form $\omega$ and $f:\C\>X$ be a $J$-holomorphic
map. We want to show $f$ is a constant, i.e. $f^*\omega=0$. There is
a lifting $\tilde f:\C\>\tilde X$ of $f$ such that $\pi_X\circ
\tilde f =f$. Let $g$ be the induced  Riemannian metric on
$(X,\omega)$ defined as in (\ref{smetric}). The induced almost complex
structure, symplectic form and metric on $\tilde X$ are denoted by $\tilde J$, $\tilde{\omega}$ and
$\tilde g$ respectively.  On the other
hand, for any tangent vectors $v, w$ on $\C$,
\be \tilde g\left( \tilde f_*v, \tilde f_*w \right)
&=&\tilde \omega\left(\tilde f_*v, \tilde{J} \tilde
f_*\left(w\right)\right)=\tilde \omega\left(\tilde f_*v, \tilde
f_*\left(J_0w\right)\right) =\left( f^*  \omega\right)
\left(v,J_0 w\right)\\&=& \omega
\left(f_{\ast}v, f_{\ast} J_0 w\right) = \omega
\left(f_{\ast}v, J  f_{\ast} w\right)=\left(f^*g\right)\left(v,  w\right),\ee
where $J_0$ is the standard complex structure of $\C$.
Hence, we have
\beq \label{4.3} \tilde f^* \tilde g =f^* g.\eeq
 Since $\pi_X^*\omega=d\eta$, \beq  f^*\omega= d
(\tilde f^*\eta) \label{po} \eeq where $\eta$ is bounded on $(\tilde
X,\tilde g)$. On the other hand, since $J$ is compatible with the
symplectic form $\omega$ and $f:\C\>X$ is $J$-holomorphic, if $f$ is
not a constant, then
$$f^{\ast}\omega=\frac{1}{2}\sqrt{-1} \mu(z) dz\wedge d\bar{z}$$ is a
K\"ahler form which may degenerate at countably many points on $\C$.
Moreover,  the  potential  $\tilde f^*\eta$ in (\ref{po}) is still
($f^{\ast}g$)-bounded on $\C$. Indeed, since $\eta$ is bounded over $(\tilde X,\tilde g)$,  for any tangent vector  $v$ on $\C$, we
have
$$\left|\left(\tilde f^*\eta\right)(v)\right|^2=\left|\eta\left(\tilde f_* v\right)\right|^2\leq C \left|\tilde f_* v\right|_{\tilde g}^2=C|v|^2_{f^{\ast }g},$$
where we use (\ref{4.3}) in the last step.

 For any bounded domain $\Omega\subset \C$,
we use $A_{\mu}(\Omega)$ and  $L_{\mu}(\partial \Omega)$ to denote
the area of $\Omega$ and the length of $\partial \Omega $ with
respect to the measure  $\mu(z) |dz|$. Then
\begin{equation}
A_{\mu}(\Omega)=\int_{\Omega} f^{\ast}\omega= \int_{\partial \Omega}
\tilde{f}^{\ast}\eta\leq C L_{\mu}(\partial \Omega)
\end{equation}
since  $\tilde{f}^{\ast}\eta$ is ($f^{\ast}g$)-bounded. Denote
$B_r=\{z\in \C:|z|<r\}$, $S_r=\{z\in \C:|z|=r\}$.   For any $r>0$,
we have \be A_{\mu}(B_r)&=&\iint_{B(r)}\mu(z)dxdy=\int_{0}^{r}
\left( \int_{S_t} \mu \right) dt\\&\geq& \int_{0}^{r}
\left(\int_{S_t} \sqrt{\mu}\right)^2(2\pi t)^{-1}
dt= \frac{1}{2\pi}\int_{0}^{r} L^2_{\mu}(S_t)\frac{1}{t}dt\\
&\geq& \frac{1}{2\pi C^2}\int_{0}^{r} A^2_{\mu}(B_t)\frac{1}{t}dt.
\ee Denote $F(r)=\int_{0}^{r} A^2_{\mu}(B_t)\frac{1}{t}dt$, then
\begin{equation}
\begin{split}
t\cdot\frac{d}{dt}F(t)& =A^2_{\mu}(B_t)\geq \frac{F(t)^2}{4\pi^2
C^4}
 \end{split}\end{equation}
 which implies
\begin{equation}
\begin{split}
\frac{d}{dt}\left(-\frac{1}{F(t)}\right) & \geq \frac{1}{4\pi^2 C^4}
\frac{1}{t}.
 \end{split}\end{equation}
 Integrating the above formula over interval $[a,b]$ with $b>a>0$, we have
 \begin{equation}
\begin{split}
\frac{1}{F(a)}\geq \frac{1}{F(b)}+\frac{1}{4\pi^2 C^4} \log
\frac{b}{a}.
 \end{split}\end{equation}
 Let $b \rightarrow \infty$, we find $F(a)=0$ for any $a>0$. This is a contradiction.
 \eproof

\noindent\emph{The proof of Theorem \ref{symp}.}  By Lemma
\ref{sectional}, we see that a compact Riemannian manifold $Y$ with
strictly negative Riemannian sectional curvature has $\tilde
d$-bounded $q^{\text{th}}$ cohomology for all  $q\geq 2 $.
 Hence, we
can apply Theorem \ref{general2} to complete the proof of Theorem
\ref{symp}.\qed

As a special case  of Theorem \ref{general2}, one can see

\bcorollary[Gromov] Let $(X,\omega)$ be a compact K\"ahler manifold.
If $X$ is K\"ahler hyperbolic, then it is Kobayashi hyperbolic.
\ecorollary

On the other hand,  the following result on Kobayashi hyperbolicity
is fundamental (e.g. \cite[Theorem~3.6.21]{Kob98}):

\btheorem Let $X$ be a compact complex manifold. If the cotangent
bundle $\Om_X^1$ is ample, then $X$ is Kobayashi hyperbolic.
\etheorem

\noindent One may wonder whether a similar result holds for K\"ahler
hyperbolicity. Unfortunately, we observe that

\bcorollary\label{coro} Let $X$ be a  complete intersection with
ample cotangent bundle $\Om_X^1$ and $\dim X\geq 2$. Then $X$ is
Kobayashi hyperbolic. However,

\bd
\item  $X$ is not K\"ahler hyperbolic.
\item $X$ can not admit a
 Riemannain metric with  non-positive sectional curvature.
\ed \ecorollary

\bproof  It is well-known that complete intersections are all simply
connected (e.g. \cite[p.225-p.227]{Sha94}).  \eproof

 \noindent For example, the
intersection of two generic hypersurfaces in $\P^4$ whose degrees
are greater than $35$ has ample cotangent bundle
(\cite[Corollary~4.13]{Bro14}). In particular, complete
intersections are not strongly negative in the sense of Siu.

\end{document}